\documentclass{amsart} 
\usepackage{amssymb, amscd}
\numberwithin{equation}{section}

\def\EE{{\mathbb E}}    
\def\FF{{\mathbb F}}  
\def\GG{{\mathbb G}}

\def\PP{{\mathbb P}}

\def\ZZ{{\mathbb Z}} 

\def\tto{\longrightarrow}

\def\Acal{{\mathcal A}}
\def\Bcal{{\mathcal B}}

\def\Gcal{{\mathcal G}}
\def\Hcal{{\mathcal H}}

\def\Lcal{{\mathcal L}}

\newcommand\riso{\mathrel{\hskip2pt\raise-2.5pt\hbox{$\widetilde{\phantom{xx}}$}
\kern-16pt\longrightarrow}}

\newcommand\proofsquare{\nobreak\hfill \hbox{%
\vrule height 5pt 
\kern-.4pt
 \vbox{%
\hrule width 5pt depth0pt height.4pt
 \kern4.6pt \hrule  }
\kern-3.75pt 
\vrule height 5pt}\kern1pt
\par}

\newtheorem{theorem}{Theorem}[section]
\newtheorem{lemma}[theorem]{Lemma}
\newtheorem{proposition}[theorem]{Proposition}
\newtheorem{corollary}[theorem]{Corollary}

\newtheorem{definition-lemma}[theorem]{Definition-Lemma}

\theoremstyle{definition}
\newtheorem{definition}[theorem]{Definition}

\theoremstyle{remark} 
\newtheorem{remark}[theorem]{Remark}

\newtheorem{question}[theorem]{Question}

\begin{document}

\title[
An invariant for varieties in positive characteristic]
{An invariant for varieties in positive characteristic
\footnote{\tt anumber.tex\today}} 
\author{G.\ van der Geer}

\address{Faculteit Wiskunde en Informatica, University of
Amsterdam, Plantage Muidergracht 24, 1018 TV Amsterdam, The Netherlands}
\email{geer@science.uva.nl}
\author{T.\ Katsura}
\address{Graduate School of Mathematical Sciences\\
 University of Tokyo \\ Komaba, Meguro-ku\\
 Tokyo\\
153-8914 Japan}
\email{tkatsura@ms.u-tokyo.ac.jp}


\begin{abstract}
We introduce an invariant of varieties in positive characteristic
which generalizes the $a$-number of an abelian variety.
We calculate it in some examples and discuss its meaning for
moduli.
\end{abstract}

\maketitle

\begin{section}{Introduction}

Varieties in positive characteristic can differ markedly from varieties
in characteristic zero. One aspect that makes this clear is Hodge theory.
If $X$ is a smooth variety in positive characteristic for which the 
Hodge-to-de Rham spectral sequence degenerates then the 
de Rham cohomology groups 
$H^m_{\rm dR}(X)$ possess two filtrations, a decreasing filtration
called the Hodge filtration $F^{\bullet}$
and an increasing filtration called the conjugate filtration 
$G_{\bullet}$.
The conjugate filtration is the analogue of the
complex conjugate of the Hodge filtration in characteristic zero.
But while in characteristic zero the Hodge filtration and the conjugate
filtration are always transversal, this is no longer the case in positive
characteristic and this gives rise to interesting new invariants of algebraic
varieties in positive characteristic. Varieties in positive characteristic
for which the Hodge filtration and the conjugate filtration are transversal
are called `ordinary' and resemble in some way smooth complex varieties,
while varieties with non-transversal filtrations resemble singular complex
varieties.

The relative position of the two filtrations on the de Rham cohomology 
is encoded by a double coset of a Weyl group and is rather complicated.
We introduce an invariant which measures the position of the first step  
of the conjugate filtration in the Hodge filtration. If $X$ is an abelian
variety then this invariant coincides with 
the $a$-number as defined by Oort.
The non-transversality is related to the cohomology of the sheaves 
$B_1\Omega^i$ introduced by Illusie.

In this paper after some preliminaries we define the $a$-number and show that
this definition coincides with the definition of the $a$-number for
abelian varieties and indicate the meaning of this number for the moduli
of abelian varieties. This is related to work of Ogus on the order
of singularity of the Hasse locus in the moduli of Calabi-Yau varieties. 
We then show how to calculate the $a$-number
for Fermat varieties using the Poincar\'e residue map. After that we
show an inequality for the $a$-number of Calabi-Yau threefolds and finish with
a discussion of the $a$-number for a family of quintic Calabi-Yau varieties.

\end{section}
\begin{section}{The Two Spectral Sequences}
Let $X$ be a smooth complete algebraic variety defined over an algebraically
closed field $k$ of characteristic $p>0$.  Let $X^{\prime}$ be the base-change
of $X$ with respect to the Frobenius homomorphism of $k$. 
The absolute Frobenius $F_{\rm abs}: X \tto X$
factors through the relative Frobenius
$ F_{\rm abs} : X {\buildrel F \over \longrightarrow} X' \tto X $.
There are two spectral sequences converging to the de Rham cohomology:
the first is the Hodge spectral sequence
$$
E_1^{ij}= H^j(X, \Omega^i) \Longrightarrow H_{\rm dR}^{i+j}(X).
$$
This spectral sequence arises from the filtration $\Omega^{\geq i}$
of the de Rham complex $\Omega^{\bullet}$.
The second spectral sequence is the so-called
 conjugate spectral sequence which comes from the 
the Leray spectral sequence for the relative Frobenius
$F: X \tto X^{\prime}$,
that yields a spectral sequence
$$
E_2^{i,j} = H^i(X^{\prime}, {\Hcal}(F_*\Omega_{X/k}^{\bullet})) \Longrightarrow
H_{dR}^{i+j}(X/k).
$$
But the Cartier operator yields an isomorphism of sheaves on $X^{\prime}$
$$
C^{-1}\colon \Omega_{X^{\prime}/k}^i {\buildrel\sim \over \tto}
{\Hcal}^i(F_*(\Omega_{X/k}^{\bullet})),
$$
so that we can rewrite this as
$$
E_2^{ij}= H^i(X^{\prime}, {\Hcal}^j(F_*(\Omega^{\bullet}_{X/k}))) =
H^i(X^{\prime},\Omega_{X^{\prime}/k}^j)
\Longrightarrow H_{\rm dR}^*(X).
$$
We assume that the Hodge-to-de Rham spectral sequence for $X$
degenerates at the $E_1$-level. This happens for example
if $p> \dim(X)$ and $X$ can be lifted to the Witt vectors $W_2(k)$, see [De-Il].
 Then the de Rham cohomology
carries a filtration, the Hodge filtration
$$
F^{\bullet}: \, H^m_{dR}=F^0\supset F^1\supset \ldots \supset F^m,
$$
with graded pieces
$$
{\rm gr}^i H^{m}_{dR}(X) = H^{m-i}(X, \Omega^i).
$$
Moreover, then also the conjugate spectral sequence degenerates at the
$E_2$-level, leading to the so-called conjugate filtration
$$
 G_{\bullet}: (0)\subset G_0\subset G_1 \subset \ldots \subset G_m=H^m_{dR}(X)
$$
with graded pieces
$$
{\rm gr}_i H^{m}_{dR}(X) = H^i(X^{\prime}, \Omega^{m-i}_{X^{\prime}}).
$$
Moreover, if $m=n:=\dim(X)$ then we have a non-degenerate pairing 
$\langle \, , \, \rangle$ on the de Rham cohomology
$$
H^n_{\rm dR}(X) \times H^n_{\rm dR}(X) \to H^{2n}_{\rm dR}(X) \cong k
$$
Note that for the conjugate spectral sequence we have
$$
E_2^{n \, 0} =H^n(X^{\prime}, {\Hcal}^0(\Omega^{\bullet}_{X^{\prime}/k}))=
 H^n(X^{\prime}, O_{X^{\prime}})$$
and since $F^*_{\rm abs}(H^n(X,O_X))=H^n(X^{\prime}, O_{X^{\prime}})$ we see that
$G_0=F_{\rm abs}^*(H^n(X,O_X))$. In particular, the composition
$$H^n(X,O_X) {\buildrel F^*\over \tto}G_0 \tto F^0/F^1 \cong H^n(X,O_X).
$$
is the Hasse-Witt map.
\par

\end{section}
\begin{section}{The relative position of two filtrations}
Given two filtrations on the de Rham cohomology it is natural to compare
them to obtain interesting information on the variety.
In general, if one has two filtrations on a finite-dimensional 
$k$-vector space $V$
$$
F^{\bullet}: V=F^0\supset F^1\supset \ldots \supset F^m, 
\quad  G_{\bullet}: G_0\subset G_1 \subset \ldots \subset G_m=V
$$
such that ${\rm rank}(G_i)={\rm rank}(F^{m-i})$
the relative position of two such filtrations is encoded by an 
element of a double coset of a Weyl group of the general linear 
group ${\Gcal}={\rm GL}(V)$.   Sometimes, the vector space $V$ is provided
with a non-degenerate pairing $\langle \, , \, \rangle $ and then we
can consider the group ${\Gcal}={\rm GSp}(V)$ instead of ${\rm GL}(V)$.

We fix a maximal torus $T$ and a Borel subgroup $B_0$ containing $T$ 
and let $W$ be the Weyl group of ${\Gcal}$. 
The (partial) flags $F^{\bullet}$ and
$G_{\bullet}$ determine a parabolic subgroups $P_{F}$ and $P_{G}$
which are conjugate under $W$. Let $W_F=W_G$ be the Weyl group of $F^{\bullet}$.
Then the relative position of $F^{\bullet}$ and $G_{\bullet}$ is 
given by an element of
$$
w(F^{\bullet}, G_{\bullet}) \in W_F\backslash W /W_G.
$$
To define it, recall that if ${\Bcal}$ denotes the variety of Borel subgroups
of ${\Gcal}$ then we have a bijection $\phi: W {\buildrel \sim
\over \longrightarrow} {\Gcal}\backslash {\Bcal} \times {\Bcal}$, given by 
$w \mapsto $ the orbit of $(B,wBw^{-1})$. We choose Borel 
subgroups $B_F \subset P_F$ and $B_G \subset P_G$ and define 
$w(F^{\bullet}, G_{\bullet})$ as the double coset containing $(B_F,B_G)$
and this is independent of the choices of $B_F$ and $B_G$.
Note that the double coset $W_F\backslash W /W_G$ is in bijection with
$P_G \backslash {\Gcal} / P_F$, cf.\ \cite{M}. 

\smallskip
\noindent
\begin{definition} Let $X$ be a smooth variety of dimension $n$ 
in characteristic $p>0$ for which
the Hodge-to-de Rham spectral sequence degenerates and let $m$ be
an integer with $0\leq m \leq 2n$. We define the pointer $w_m(X)$
of $X$ as the element of $W_F\backslash W /W_G$ associated to the
two filtrations $\{ F^{\bullet}\}$ and $\{ G_{\bullet} \}$ on $H^m_{\rm dR}(X)$.
\end{definition}

\par
We can refine this considerably by choosing full filtrations refining the
filtrations $F^{\bullet}$ and $G_{\bullet}$ 
$$
\begin{matrix}
V=\Phi^0\supset \Phi^1\supset \ldots \supset \Phi^{\dim(V)},&\\
\Gamma_0 \subset \Gamma_1\subset \ldots \subset \Gamma_{\dim(V)}=V&\\
\end{matrix}
$$
that are compatible with the action of Frobenius in the following sense.
Using crystalline cohomology we know that $F/p^i$ acts on $F^i/F^{i+1}$
and the induced filtration of $\Phi^{\bullet}$ should be stable under $F/p^i$.
Moreover, the filtration $\Phi^{\bullet}$ induces via
$G_i/G_{i-1} \cong F^{i-1}/F^i$ the filtration $\Gamma_{\bullet}$.
Comparing the two filtrations gives the finer invariant. If we carry this out
for principally polarized abelian varieties we retrieve the Ekedahl-Oort
type of the abelian variety, cf.\ \cite{O2} where this invariant was introduced
in a different way. We refer to  \cite{G} and to \cite{E-G}, where this
will be worked out in detail for abelian varieties, K3-surfaces and 
Calabi-Yau's and leads to stratifications on moduli spaces (cf.\ also
\cite{M}, \cite{G-K}).
\end{section}
\begin{section}{The $a$-number of a variety in positive characteristic}
We are now interested in the case of the cohomology group
$H^m_{dR}(X)$ with $m=n=\dim(X)$. We shall assume throughout that the
Hodge-to-de Rham spectral sequence degenerates at the $E_1$-level.
Since the pointer $w_n(X)$ is a complicated
invariant we look first at the position of $G_0$, the image of $H^n(X,O_X)$
under the Frobenius operator in the Hodge filtration $F^{\bullet}$
on $H^n_{dR}(X)$.  We distill an invariant in the following way.
\bigskip
\noindent
\begin{definition} The $a$-number of $X$ is the maximum
step in the Hodge fitration that contains $F^*_{\rm abs}(H^n(X,O_X))$:
$$
a(X):=\max \{ j : F^jH^n_{dR}(X) \, \hbox{\rm contains} \, 
F_{\rm abs}^*(H^n(X,O_X)) \},
$$
or equivalently,
$$
a(X)=\max \{ j: F^j \supseteq G_0\}.
$$
\end{definition}
Note that  $0 \leq a(X) \leq n=\dim(X)$ and that $a(X)>0$ if and only 
if the Hasse-Witt map vanishes.
Of course, by going from $w_n(X)$ to $a(X)$ we loose a lot of information.
If we wish to retain more information we can refine the $a$-number slightly
by taking the $a$-vector  $(a_0,a_1,\ldots,a_m)(X)$ with 
$$
a_j(X):=\dim G_0 \cap F^j \quad j=0,1,\ldots,m.
$$ 
Sometimes, if $H^n(X,O_X)$ vanishes, we can still define an analogue
of the $a$-number. Consider for example the case of a smooth cubic hypersurface
in $\PP^8$. The Hodge numbers $h^{ij}$ with $i+j=7$ are
$$
h^{0 , 7}=h^{7 , 0}=0, \, h^{1 , 6}=h^{6 , 1} = 0, \, h^{2 , 5}=
h^{5 , 2}=1, \, h^{3 , 4}=h^{4 , 3}= 360.
$$
and we can define
$$
a^{\prime}(X)= \max\{ j: F^j \supseteq G_2\}.
$$
\begin{question}
i) Can we give geometric interpretations of these numbers?
ii) Is a variety with $a(X)=\dim(X)$ rigid ? ; i.e., are there deformations
that preserve the $a$-number?
\end{question}
\end{section}
\begin{section}{ Closed differential forms}
\bigskip
Let $X$ be a smooth projective variety of dimension $n$ 
in characteristic $p>0$
and let $F_{\rm abs}: X \tto X$ be the absolute Frobenius.
We consider the direct image ${F_{\rm abs}}_*O_X$. As a sheaf of abelian groups
it is just $O_X$, but its $O_X$-module structure is different:
$f \circ g = f^p\, g$. One way to define the sheaf $B_1\Omega^{i}_X$ is via
the short exact sequence
$$
0 \to \Omega_X^{i-1} \tto {F_{\rm abs}}_* \Omega_X^{i-1}
{\buildrel d \over \longrightarrow}
B_1\Omega_X^{i} \to 0.
$$
Here we view $B_1\Omega_X^{i}$ as an $O_X$-module. It can also be viewed as
a locally free subsheaf of $F_*\Omega_X^i$ on $X^{\prime}$. Here $F$ denotes
the relative Frobenius.
Moreover, we set $B_0\Omega_X^i=0$. Inductively we define
$$
B_{j+1}\Omega_X^i=C^{-1}(B_j\Omega_X^i),
$$
where $C: Z_1\Omega_X^i \to \Omega_X^i$ is the Cartier operator.
Similarly, we define
$$
Z_0\Omega_X^i= \Omega_X^i, \, Z_1\Omega_X^i=
\Omega_{X,\hbox{\rm d-closed}}^i
$$
and inductively
$$
Z_{j+1}\Omega_X^{i}=C^{-1}(Z_j\Omega_X^i)\qquad \hbox{\rm for $j\geq 1$}.
$$
We can view these as locally free subsheaves of
$F^j_* \Omega_X^i$ on $X^{(p^j)}$, the base change of $X$ under the $j$-th
power of relative Frobenius.
The inverse Cartier operator gives rise to an isomorphism
$$
C^{-j}: \Omega_{X^{(p^j)}}^i {\buildrel \cong \over \longrightarrow}
Z_j\Omega_X^i/B_j \Omega_X^i.
$$
We have a perfect pairing
$$
{F_{\rm abs}}_*\Omega_X^i \otimes {F_{\rm abs}}_*\Omega_X^{n-i} \tto \Omega_X^n,
\qquad (\omega_1 , \omega_2) \mapsto C(\omega_1 \wedge \omega_2).
$$
On the other hand we have an exact sequence
$$
0 \to \Omega_X^i \tto {F_{\rm abs}}_*\Omega_X^i
{\buildrel d \over \longrightarrow}  B_1\Omega_X^{i+1} \to 0.
$$
and
$$
0\to B_1\Omega_X^{n-i} \tto {F_{\rm abs}}_*\Omega_X^{n-i} {\buildrel C \over
\longrightarrow} \Omega_X^{n-i} \to 0.
$$
This induces a perfect pairing
$$
B_1\Omega_X^{i+1} \otimes B_1\Omega^{n-i}_X \to \Omega_X^n .
$$
\begin{definition} (Illusie, Raynaud)
We call the variety $X$ {\sl ordinary} if the cohomology groups
$ H^i(X, B_1\Omega^j_X)$ vanish for $j\geq 1$ and all $i$.
\end{definition}
This implies that all global forms are closed, just as in characteristic $0$.
\begin{lemma} If $X$ is ordinary then $a(X)=0$. 
\end{lemma}
\begin{proof} Consider the exact sequence
$$
0\to O_X {\buildrel F \over \tto} O_X {\buildrel d \over \tto} B_1\Omega^1_X
\to 0.
$$
In cohomology this gives 
$$
H^{n-1}(X, B_1\Omega_X^1) \to 
H^n(X,O_X) \to H^n(X,O_X) \to H^n(X, B_1\Omega_X^1)
$$
from which it follows that the Hasse-Witt map has trivial kernel and thus that
$G_0 \not\subseteq F^1(H^n_{\rm dR}(X))$.
\end{proof}
\end{section}
\begin{section}{Abelian Varieties}
Let $X$ be an abelian variety of dimension $g$ over an algebraically 
closed field $k$. 
We denote by $X[p]$ the kernel of multiplication
by $p$. It is a group scheme of order $p^{2g}$.
The classical $a$-number of $X$ (cf.\ \cite{O1}) is defined by
$$
a(X):= \dim_k {\rm Hom}_k(\alpha_p, X).
$$
Here $\alpha_p$ is the group scheme of order $p$ 
that is the kernel of Frobenius acting on the additive group $\GG_a$.
Note that ${\rm Hom}_k(\alpha_p,X)$, the space of group scheme
homomorphisms of $\alpha_p$ to $X$, is in a natural way a vector
space over $k$. If we let $A(X)$ be the maximal subgroup scheme of $X[p]$
annihilated by the operators $F$ (Frobenius) and $V$ (Verschiebung), 
then $A(X)$ is the union of the images of all group scheme 
homomorphisms $\alpha_p \to X$ and we have
$$
a(X)= \log_p {\rm ord} A(X).
$$
We have $0 \leq a(X) \leq g$.

It is well-known that the contravariant Dieudonn\'e module $D(X)$
of $X[p]$ can be identified with the first de Rham cohomology 
$H^1_{\rm dR}(X)$ of $X$. The $k$-vector space $H^1_{\rm dR}(X)$ then carries 
two operators $F$ and $V$ and the Dieudonn\'e module of $A(X)$
coincides with the intersection $\ker(V) \cap \ker(F)$ on $D(X)$.
If we identify the kernel of $F$ with 
$H^0(X, \Omega_X^1) \subset H^1_{\rm dR}(X)$
then the Dieudonn\'e module of $A(X)$ may be identified with the
kernel of $V$ acting on $H^0(X, \Omega_X^1)$. We thus find
$$
D(A(X)) \cong \ker(V: H^0(X,\Omega_X^1)\to H^0(X,\Omega_X^1)) 
\cong H^0(X, B_1\Omega_X^1).
$$
and we have
$$
a(X)= \dim_k H^0(X,B_1\Omega_X^1).
$$
We also have the following relation:
$$
F(H^1(X,O_X)) \cap H^0(X, \Omega_X^1)= H^0(X,B_1 \Omega_X^1). \eqno(1)
$$
This follows from the fact that ${\rm Im}(F)={\rm Ker} (V)$.
\smallskip
We now show that for abelian varieties the classical 
$a$-number and our $a$-number coincide.

\begin{proposition} If $X$ is an abelian variety of dimension $g$
then the two notions of $a$-number coincide: if 
$F^{\bullet}$ and $G_{\bullet}$ are the Hodge and conjugate
filtration on $H^g_{\rm dR}(X)$ 
we have
$$
\dim_k {\rm Hom}_k(\alpha_p,X) = \max \{ j: G_0 \subset F^j\}.
$$
\end{proposition}
\begin{proof} Recall that $H_{\rm dR}^g(X) = \wedge^g H^1_{\rm dR}(X)$ and
if we write $H^1_{\rm dR}(X)=V_1\oplus V_2$ with $V_1=H^0(X, \Omega_X^1)$
and $V_2$ a complementary subspace, then
the Hodge filtration on $H^g_{\rm dR}(X)$ is 
$F^r= \sum_{j=r}^g \wedge^j V_1 \otimes \wedge^{g-j} V_2$.
We have  $F(H^g(X,O_X))= F(\wedge^g H^1(X,O_X))=
\wedge^g F(H^1(X,O_X))$. If we write $F(H^1(X,O_X))= A\oplus B$
with $A=H^0(X,B_1\Omega_X^1)$ and $B$ a complementary space, then
$\wedge^g (A\oplus B)= \wedge^a A \otimes \wedge^{g-a} B$ with $a=\dim (A)$.
From this and (1) it is clear that $F(H^g(X,O_X))$ lies in
$F^a$, but not in $F^{a+1}$. \end{proof}
We now show that the $a$-number has some meaning for the geometry of
moduli spaces.

Let $T(a)$ be the locus inside the moduli space ${\Acal}_g$ 
of principally polarized abelian varieties with $a$-number $\geq a$.
Here ${\Acal}_g$ is viewed as an algebraic stack, or we should
add a sufficient level structure to our abelian varieties. 
Over ${\Acal}_g$ the de Rham bundle ${\Hcal}_{\rm dR}^1$
posseses two subbundles of rank $g$, the Hodge bundle $\EE$ and 
the kernel $\FF$ of $F$. Then $T(a)$ is defined as the degeneracy locus
where $\EE \cap \FF$ has rank at least $a$. 
It is known that $\dim T(a) = g(g+1)/2- a(a+1)/2$, cf.\cite{G}, 
\cite{E-G}. 
\smallskip
\begin{proposition} The locus $T(a)$ in ${\Acal}_g$ is smooth outside $T(a+1)$ 
and the normal space to $T(a)$ at $[X]$  with $a(X)=a$
can be identified with ${\rm Sym}^2(H^0(X,B_1\Omega_X^1))$.
\end{proposition}
\begin{proof} 
An infinitesimal deformation of the principally polarized abelian
variety $X$ is given by a symmetric $g\times g$-matrix $T=(t_{ij})$
which can be interpreted as a symmetric endomorphism of $H^0(X,\Omega^1_X)$,
cf.\ \cite{Ko}. 
This deformation preserves the $a$-number of $X$ if it keeps the kernel
$H^0(X,B_1\Omega_X^1)$ of $V$ acting on $H^0(X,\Omega_X^1)$, that is,
$T$ is a symmetric endomorphism of this subspace. The principal
polarization identifies this subspace with its dual, and 
this gives us the result.
\end{proof}

\begin{proposition}
Let $X$ be an abelian variety with $a(X)=g$.
Then the multiplicity of the point $[X]$ on $T(1)$ is $g$. \end{proposition}
\begin{proof} If $X$ is a principally polarized 
abelian variety with $a(X)=g$ we 
choose a base $\omega_1,\ldots,\omega_g$ of $H^0(X,\Omega_X^1)$ and
complete it to a basis of $H^1_{\rm dR}(X)$ as
$\eta_1,\ldots,\eta_g$ such that $\langle \omega_i,\eta_j\rangle =\delta_{ij}$.
Then for a deformation with parameter a symmetric $g\times g$-matrix $(t_{ij})$
we have for $i=1,\ldots,g$
$$
\begin{aligned}
F\eta_i &= \omega_1 + \sum t_{ij} \eta_j\\
\omega_i &= V (\sum t_{ij}\eta_i).\\
\end{aligned}
$$
Then we get
$$
F(\eta_1\wedge \cdots \wedge \eta_g)= \omega_1\wedge\cdots \wedge \omega_g +
\cdots + \det(t_{ij}) \eta_1\wedge \cdots \wedge \eta_g.
$$
But the equation of $T(1)$ is 
$$
\langle F(\eta_1\wedge \cdots \wedge \eta_g),
\eta_1\wedge \cdots \wedge \eta_g\rangle = \det(t_{ij}).
$$
\end{proof}
\end{section}
\begin{section}{Cohomology of Hypersurfaces}
Let $X$ be an irreducible smooth projective hypersurface of degree 
$d$ in $\PP^{n+1}$ given by an equation $f=0$. Then the primitive 
cohomology $H^n_{\rm dR}(X)_0$ in the middle dimension can be described 
by the Poincar\'e residue map as follows.  Consider the rational
 differential forms on projective space
$$
\Omega= \sum_{i=0}^{n+1} (-1)^i x_i dx_0\wedge \ldots \wedge d\hat{x}_i \wedge
\ldots \wedge dx_{n+1} 
\quad\hbox{\rm and} \quad 
\Omega^*= \frac{\Omega}{x_0\ldots x_{n+1}}.
$$
Note that $\Omega^*$ is a logarithmic form: in affine coordinates we have
$\Omega^*= {du_1}/{u_1}\wedge \ldots \wedge {du_{n+1}}/{u_{n+1}}$.

We let ${\Lcal}$ be the $k$-vector space generated by the monomials
$x^w= x_0^{w_0}\ldots x_{n+1}^{w_{n+1}}$ with $\sum_{i=0}^{n+1} w_i 
\equiv 0 (\bmod \, d)$ and we let ${\Lcal}^{\prime} \subset {\Lcal}$ be 
the subspace generated by the $x^w$ with all $w_i \geq 1$, i.e., the 
elements in ${\Lcal}^{\prime}$ are divisible by $x_0\ldots x_{n+1}$. 
For an element $x^w \in {\Lcal}^{\prime}$ we set $\gamma(w)= \sum w_i/d$.

On ${\Lcal}$ we have operators $D_i$ defined by
$$
D_i(x^w)= 
x_i \frac{\partial x^w}{\partial x_i} +x_i \frac{\partial f}{\partial x_i}x^w=
w_ix^w + x_i \frac{\partial f}{\partial x_i}x^w
. \eqno(2)
$$
There is a natural map
$$
\phi:{\Lcal}^{\prime} \longrightarrow H_{\rm dR}^{n+1}(\PP^{n+1}\backslash X),
\qquad x^w \mapsto (-1)^{\gamma -1} (\gamma -1)! \frac{x^w}{f^{\gamma}}
\Omega^*.
$$
Now we use the Poincar\'e residue map
$$
 {\rm Res}: H_{\rm dR}^{n+1}(\PP^{n+1}\backslash X) \longrightarrow
H^n_{\rm dR}(X)
$$
and observe that the composition $\rho:= {\rm Res} \cdot \phi$ factors through 
quotient
$$
W^{\prime}={\Lcal}^{\prime}/{\Lcal}^{\prime}\cap \sum_{i=0}^{n+1} D_i {\Lcal}. 
$$
and we get an isomorphism $W^{\prime} \cong H^n_{\rm dR}(X)_0$
which is compatible with the pole order filtration on the source 
and the Hodge filtration on the target, cf.\ \cite{K}, cf.\ also
\cite{C-G},\cite{St}.
The image of $H^n(X,O_X)$ under the
action of absolute Frobenius on $H^n_{\rm dR}(X)$ is up to a multiplicative
constant given by raising $x^w/f^{\gamma}$ to the $p$-th power
and then reducing the pole order by making use of the relations (2).
\end{section}
\begin{section}{Fermat Varieties}
We give examples of surfaces with $a$-number equal to $0,1$ or $2$.
\begin{theorem} Let $p$ be a prime different from $5$.
Let $X$ be the Fermat surface in characteristic $p$
 defined by $\sum_{i=0}^3 x_i^5=0$ in $\PP^3$.
Then the $a$-number of $X$ satisfies
$$
a(X)= \begin{cases}
0 & p\equiv 1 (\bmod \, 5)\\
1 &  p\equiv 2,\, 3 (\bmod \, 5)\\
2 &  p\equiv 4 (\bmod\,  5).\\
\end{cases}
$$
\end{theorem}
\begin{proof}
We use the description of the (primitive) cohomology  $H^2(X,O_X)$ 
with the Poincar\'e residue map
$$
P:V \tto H_{\rm dR}^2(X), \qquad A \mapsto {\rm Res}(
\frac{A\, x_0x_1x_2x_3}{f^3} \Omega^*)
$$
with $V$ the vector space of homogeneous polynomials of degree $11$. 
A basis of the $4$-dimensional space 
$H^2(X,O_X)$ is given by $\alpha_i=P((x_0x_1x_2x_3)^3
/x_i)$ for $i=0,\ldots,3$. Since $\Omega^*$ is a logarithmic form
the image of $\alpha_0$ under Frobenius is given by the Poincar\'e residue
of
$$ 
\frac{(x_0x_1x_2x_3)^{3p} (x_1x_2x_3)^p }{f^{3p}} \Omega^*
$$
and using the relations
$-5 x_i^5 x^w= w_ix^w \quad \hbox{\rm for any monomial $x^w$}$ obtained from (2)
this is equivalent to an expression ${\rm res} (g \Omega^*)$
with $g$ given by
$$
\begin{matrix}
(x_0x_1x_2x_3)^3 x_1x_2x_3 /f^3 & p\equiv 1 (\bmod \, 5)\\
(x_0x_1x_2x_3) (x_1x_2x_3)^2 /f^2 & p\equiv 2 (\bmod \, 5)\\
x_0^4 (x_1x_2x_3)^2/f^2 & p\equiv 3 (\bmod \, 5)\\
x_0^2 x_1x_2x_3 /f & p\equiv 4 (\bmod \, 5)\\
\end{matrix}
$$
and similarly for the other $\alpha_i$.
This pole order implies that $a(X)$ is as indicated.

\end{proof}

\begin{theorem}
Let $X$ be the Calabi-Yau Fermat variety in $\PP^r$ given by the
equation $x_0^{r+1}+\ldots+x_r^{r+1}=0$. Then the $a$-number of $X$
is the natural number $a$ with $0\leq a \leq r-1$ and
$a\equiv p-1 (\bmod \, r+1)$.
\end{theorem}
\begin{proof}
Under the Poincar\'e residue map $(x_0\ldots x_r)^r\Omega^*/f^r$
maps to a generator of $H^{r-1}(X,O_X)$. The image under Frobenius
is $(x_0\ldots x_r)^{pr}/f^{pr}\times \Omega^*$ and using again
the relations (2) we can reduce the pole order $pr$ modulo $r+1$.
But $pr \equiv -p (\bmod \, r+1)$ and from this the result easily
follows.
\end{proof}
There exist varieties of arbitrary dimension with maximal $a$-number
namely abelian varieties which are products of supersingular
elliptic curves.
Besides these abelian varieties the following Calabi-Yau varieties
give examples of such varieties.
\begin{corollary}
Let $X$ be the $(p-1)$-dimensional Fermat variety defined by
$$
x_0^{p+1}+x_1^{p+1}+\cdots + x_p^{p+1} =0
$$
in projective space of dimension $p$. Then the $a$-number of the
Calabi-Yau variety $X$ is equal to $p-1=\dim(X)$.
\end{corollary}

Let $X$ be the $7$-dimensional cubic in $\PP^8$ 
defined by $\sum_{i=0}^8 x_i^3=0$ in $\PP^8$.
Then $H^7(X,O_X)$ and $H^6(X,\Omega^1_X)$ vanish, while $H^5(X,\Omega_X^2)$
is $1$-dimensional. Let $a^{\prime}(X)$ be the maximum $j$ such that
the Hodge step $F^j$ contains $G_2$ (which can be seen
as an image of $H^5(X,\Omega^2)$ 
under Frobenius divided by $p^2$ in $H^2_{\rm dR}(X)$, cf.\ \cite{Katz1}).

\begin{theorem}
Let $X$ be the $7$-dimensional cubic in $\PP^8$ 
defined by $\sum_{i=0}^8 x_i^3=0$.
Then we have:
$$
a^{\prime}(X)= \begin{cases}
2 & p\equiv 1 (\bmod \, 3)\\
5 &  p\equiv 2(\bmod \, 3).\\
\end{cases}
$$
\end{theorem}
\begin{proof} Using again the description of primitive cohomology with
the Poincar\'e residue map we consider the form
$$
\frac{(x_0\ldots x_8)^2}{f^6} \, \Omega^*
$$
and then see that the image under Frobenius is equivalent to
the Poincar\'e residue of $g\, \Omega^*$
with $g$ given by
$(x_0\ldots x_8)^2/f^6$ 
if $p\equiv 1 (\bmod \, 3)$ 
and $(x_0\ldots x_8)/f^3$ 
if $p\equiv 2 (\bmod \, 3)$.
\end{proof}

\end{section}

\begin{section}{Calabi-Yau Varieties}
Let $X$ be an $n$-dimensional Calabi-Yau variety. Such a variety has a
height $h(X)$. There are several possible definitions of this, e.g. with
formal groups. Here we take as definition
$$
h-1:= \max \{ \dim H^{n-1}(X,B_i\Omega_X^1) : i \in \ZZ_{\geq 0} \}.
$$
One can use the methods of \cite{G-K} to 
see that this definition agrees with the definition using the
formal group $\Phi^n$ associated to $X$.
The natural inclusion $B_i\Omega_X^1 \subseteq \Omega^1_X$ induces a natural
map $\beta_i:H^{n-1}(X,B_i\Omega^1_X) \to H^{n-1}(X,\Omega_X^1)$.
\begin{lemma} If the natural map $\beta_i: 
H^{n-1}(X,B_i\Omega^1) \to H^{n-1}(X,\Omega^1)$
is not injective, then $\lim_{i\to \infty} \dim H^{n-1}(X,B_i\Omega_X^1)
=\infty$.
\end{lemma}
\begin{proof} The argument is similar to the cases of K3 surfaces, cf.\ 
\cite{G-K}.We give the argument for $i=1$.
Suppose that $\beta_1$ has a non-trivial 
kernel represented by a cocycle $\{ df_{\alpha_1,\ldots, \alpha_n} \}$.  
Then we have a relation
$$
df_{\alpha_1,\ldots, \alpha_n} 
= \sum (-1)^j \omega_{\alpha_1,\ldots,\hat{\alpha}_j,\ldots ,\alpha_n}.
$$
Since the Cartier map $H^0(U,\Omega^1_{U,\rm closed}) \to H^0(U,\Omega_U^1)$
is surjective on affine sets $U$ we can find closed forms
$\tilde{\omega}_{\alpha_1,\ldots,,\hat{\alpha}_j,\ldots , \alpha_n}$ and regular functions 
$g_{\alpha_1,\ldots, \alpha_n}$ on $\cap U_{\alpha_i}$ and obtain 
a new relation
$$
f_{\alpha_1,\ldots, \alpha_n}^{p-1}df_{\alpha_1,\ldots, \alpha_n}+
dg_{\alpha_1,\ldots, \alpha_n} =\sum (-1)^j 
\tilde{\omega}_{\alpha_1,\ldots,,\hat{\alpha}_j,\ldots , \alpha_n}.
$$
One now checks directly that the cocycle at the 
left hand side represents an element
of $H^{n-1}(X,B_2\Omega_X^1)$ which does not lie in the image of 
$H^{n-1}(X,B_1\Omega_X^1)$. One also checks that this element is non-zero.
The argument for other $i$ is similar.
\end{proof}

\begin{proposition}
If $h\neq \infty$ then $1 \leq h \leq h^{1,n-1}+1$. 
\end{proposition}
\begin{proof} If $h\neq \infty$ all the maps $\beta_i$ are injective. 
\end{proof}
\begin{remark} We have $h^1(\Theta)=h^1(\Omega^{n-1})$, so we a priori could
expect a stratification of the moduli with $h$ steps; but the
example of K3 surfaces shows that $h$ cannot assume all values
between $1$ and $h^{1,1}+1$. Therefore one should find bounds for~$h$.
\end{remark}
\begin{proposition}
Let $X$ be a Calabi-Yau variety such that the Hodge-to-de Rham spectral
sequence degenerates. Then we have either
$h(X)=1$ and $a(X)=0$, or
$2\leq h(X)<\infty$ and $a(X)=1$, or
$h(X)=\infty$ and $a(X)\geq 1$.
\end{proposition}
\begin{proof}
We have $h(X)=1 $ if and only if 
the map $F_{\rm abs}: H^n(X,O_X) \to F^0/F^1 =H^n(X,O_X)$ is non-zero 
and this is equivalent to $F_{\rm abs}(H^n(X,O_X)) \not\subset F^1$. So 
we have $h(X)\geq 2 $ if and only if 
the map  $F$ maps $H^n(X,O_X)$ to $F^1$ and then we have the projection
$H^n(X,O_X) \tto F^1\tto F^1/F^2 = H^{n-1}(X,\Omega_X^1)$. Now recall that we
have an isomorphism $d: O_X/FO_X \cong B_1\Omega_X^1$ and from the exact sequence
$$
0\to O_X {\buildrel F \over \longrightarrow} O_X \tto O_X/FO_X \to 0
$$
we thus get $H^{n-1}(X,B_1\Omega_X^1)\cong H^n(X,O_X)$.
We claim that the image of $H^n(X,O_X)$ under $F$ in $F^1/F^2$
is the image of $H^{n-1}(X,B_1\Omega_X^1)$ in $H^{n-1}(X,\Omega_X^{1})$.
This can be checked by a direct computation: if $\{ f_I\}$ is a
Cech $n$-cocycle representing a generator $g$ of $H^n(X,O_X)$ 
then $F(g)$ is represented by $\{ f^p_I\}$ and this is cohomologous
to $0$ in $H^n(X,O_X)$. So there exists an $(n-1)$-cocycle $h_J$
such that $\delta(h)=f^p$. The image of $g$ in $F^1/F^2$ is represented 
by $dh_J$.
Now use the proof of the preceding Proposition. \end{proof}
\begin{corollary}
Let $X$ be the Fermat Calabi-Yau variety of degree $r+1$ in $\PP^r$.
Assume that the characteristic $p$ does not divide $r+1$ and
satisfies $p \not\equiv 2 (\bmod \, r+1)$.
Then the height $h(X)$ is equal to $1$ or $\infty$.
Moreover, $h(X)=1$ if and only if $p\equiv 1 (\bmod \, r+1)$.
\end{corollary}
\begin{remark}
The corollary holds also without the condition that $p\not \equiv 2
(\bmod \, r+1)$, but for the proof we then need to use Jacobi sums.
\end{remark}
\end{section}
\begin{section}{Relation with a Result of Ogus}
In \cite{Og} Ogus proved a result on the order of vanishing of the Hasse-invariant
of a family of Calabi-Yau varieties that is closely related to
the $a$-number introduced here. Let $f: X\to S$ be a family of Calabi-Yau varieties 
such that for each fibre $X_s$ the Hodge-to-de Rham spectral sequence
degenerates and such that the Kodaira-Spencer mapping $T_{S/k} \to
R^1f_*(T_{X/S})$ is surjective. Then Ogus proves that the Hasse-invariant
vanishes to order $i$ at $s$ if $G_0(X_s) \subset F^i(X_s)$. Moreover, under
an additional natural  assumption on $X/S$  
the order of vanishing is exactly equal to the
$a$-number of $X_s$.

As an example we consider the family of Calabi-Yau hypersurfaces
 in $\PP^4$ given by
$$
X_{\alpha} : \qquad \sum x_i^5 - 5\alpha x_1\cdots x_5=0.
$$
We have $h^{1,1}=1$, $h^{2,1}=101$ and $b_3=204$. 
The Hasse-Witt invariant can be calculated (cf.\ \cite{Katz1}, \S 2.3)
$$
H(\alpha)= \sum_{m=0}^{[(p-1)/5]} \frac{(5m!) }{ (m!)^5} \alpha^{p-1-5m}.
$$
This is a polynomial of degree $p-1$ in $\alpha$. For $\alpha=0$ we find
by using Theorem 8.2 and the formula for $H(\alpha)$ 
\begin{proposition}
Let $p$ be a prime $\neq 5$. Then the $a$-number of $X_0$ is determined by
$$
a(X_0)={\rm ord}_0 H(\alpha)\equiv p-1 (\bmod \, 5).
$$
\end{proposition}

We have an action of the symmetric group 
$S_5$  by permutation of the coordinates
and of $\mu_5^3$ given by generators
$$
(x_1,x_2,x_3,x_4,x_5)
\begin{matrix}
& {\buildrel g_1 \over \longrightarrow} (x_1,\zeta x_2, x_3, x_4, \zeta^4 x_5)
\\
& {\buildrel g_2 \over \longrightarrow} (x_1,x_2, \zeta x_3, x_4, \zeta^4 x_5)
\\
& {\buildrel g_3 \over \longrightarrow} (x_1, x_2, x_3, \zeta x_4, \zeta^4 x_5)
\\
\end{matrix}
$$
Then $Y_{\alpha}= X_{\alpha}/\mu_5^3$ is a mirror family with
$h^{1,1}=101$, $h^{2,1}=1$ and $b_3=4$. We see a 4-step filtration
$$
F^3 \subset F^2 \subset F^1 \subset H_{\rm dR}^3.
$$
It is not difficult to see, with the help of a non-trivial
trace map $H^i(X_{\alpha},\Omega_{X_{\alpha}}^j) \to
H^i(Y_{\alpha}, \Omega_{Y_{\alpha}}^j)$,
that the $a$-number of $Y_{\alpha}$ equals that of $X_{\alpha}$. 
\end{section}
\begin{section}{Acknowledgement}
This research was made possible by a JSPS-NWO grant. The second
author would like to thank the University of Amsterdam and  
the first author would like to thank Prof.\ Ueno for inviting him
to Kyoto, where this paper was finished.
We also thank Ben Moonen for his comments. 
\end{section}

\end{document}